\documentclass[10pt,twocolumn,twoside]{IEEEtran}
\usepackage{bm}   % for bold
\usepackage{amsmath}
\allowdisplaybreaks[4]
\usepackage{amssymb}
\usepackage{indentfirst}  % for indent
\usepackage{graphicx}  % for pics
\usepackage{subfigure}
\usepackage{multirow} % for tabs
\usepackage{booktabs} % for tabs
\usepackage{cite}  % for citation and reference
\usepackage{threeparttable}   % for annotation in tabs
\usepackage{algorithm}% http://ctan.org/pkg/algorithm
\usepackage{algpseudocode}% http://ctan.org/pkg/algorithmicx
\usepackage{arydshln}
\usepackage{color}
\usepackage[american]{babel}
\usepackage{microtype}
\usepackage[numbers,sort&compress]{natbib}  % to include multi-citation in one []
\usepackage{url}

\usepackage{cases}
\usepackage{enumitem}
\usepackage{mathtools}
\usepackage{epstopdf}
\newtheorem{theorem}{Theorem}

\newtheorem{lemma}{Lemma}
\newtheorem{definition}{Definition}

\newtheorem{assumption}{Assumption}
\newtheorem{conjecture}{Conjecture}

\def\mV{\mathcal{V}}

\def\mbR{\mathbb{R}}
\def\mbC{\mathbb{C}}

\def\mI{\mathcal{I}}

\def\upj{\textup{j}}
\def\upeq{\textup{eq}}
\def\Jalg{\bm{J}_{\textup{alg}}}
\def\Ivs{\mI_{\textup{vs}}}
\def\mot{\textup{mot}}
\def\stat{\textup{stat}}
\def\bus{\textup{bus}}
\def\gen{\textup{gen}}

\newcommand{\diag}[1]{\textup{diag}\{#1\}}

\newcommand{\lbar}[1]{\overline{#1}}
\begin{document}
\title{Impasse Surface of Differential-Algebraic Power System Models: An Interpretation Based on Admittance Matrices}

\author{Yue~Song,~\IEEEmembership{Member,~IEEE,}
        David~J.~Hill,~\IEEEmembership{Life Fellow,~IEEE,}
        Tao~Liu,~\IEEEmembership{Member,~IEEE}\\
        and~Xinran~Zhang,~\IEEEmembership{Senior Member,~IEEE}% <-this % stops a space

%\vspace{-16pt}

%\thanks{This work was supported by the Hong Kong RGC General Research Fund under Project No. 17207918.}
\thanks{Y. Song is with the State Key Laboratory of Intelligent Autonomous Systems and Frontiers Science Center for Intelligent Autonomous Systems, Tongji University, Shanghai, China (e-mail: songy31@163.com).}
\thanks{T. Liu is with the Department of Electrical and Electronic Engineering, The University of Hong Kong, Hong Kong (e-mail: taoliu@eee.hku.hk).}
\thanks{D. J. Hill is with the Department of Electrical and Electronic Engineering, The University of Hong Kong, Hong Kong, and also with the Department of Electrical and Computer Systems Engineering, Monash University, Melbourne, Australia (e-mail: dhill@eee.hku.hk).}
\thanks{X. Zhang is with the School of Automation Science and Electrical Engineering, Beihang University, Beijing, China (e-mail: zhangxr07@buaa.edu.cn).}
}

\markboth{ACCEPTED BY IEEE TAC} %, ~Vol.~, No.~, ~2015}
%\markboth{network-based small-disturbance stability-Yue Song}
{Song \MakeLowercase{\textit{et al.}}: Impasse Surface of Differential-Algebraic Power System Models}

\maketitle

\begin{abstract}
  The impasse surface is an important concept in the differential-algebraic equation (DAE) model of power systems, which is associated with short-term voltage collapse.
  This paper establishes a necessary condition for a system trajectory hitting the impasse surface.
  The condition is in terms of admittance matrices regarding the power network, generators and loads, which specifies the pattern of interaction between those system components that can induce voltage collapse.
  It applies to generic DAE models featuring high-order synchronous generators, static loads, induction motor loads and lossy power networks.
  We also identify a class of static load parameters that prevent power systems from hitting the impasse surface; this proves a conjecture made by Hiskens that has been unsolved for decades.
  Moreover, the obtained results lead to an early indicator of voltage collapse and a novel viewpoint that inductive compensation to the power network has a positive effect on preventing short-term voltage collapse, which are verified via numerical simulations.
\end{abstract}


% \fi
\IEEEpeerreviewmaketitle

%IEEEhowto:kopka

\section{Introduction}\label{secintro}
The dynamical behaviors of electric power systems, especially those considering short-term voltage dynamics, are commonly described by a group of differential-algebraic equations (DAEs).
In the DAE model, the differential equations refer to the dynamics of synchronous generators and induction motors, while the algebraic equations refer to the power flow equations describing the balance between power transfer and load consumption~\cite{van2007voltage}.
Short-term voltage stability is a major concern in the study of power systems described by DAE models.
Its timescale is in the order of several seconds involving the dynamics of fast acting load components \cite{kundur2004definition}.
It is reported that restorative loads (e.g., induction motors) are a driving factor for short-term voltage collapse \cite{potamianakis2006short} that may cause severe damage to power systems.

\indent
Apart from the load-side viewpoint, short-term voltage collapse is closely connected to a system-wide property of DAE models, namely the impasse surface.
An impasse surface refers to the hypersurface where the algebraic Jacobian (i.e., the Jacobian matrix of the algebraic equations with respect to algebraic variables) becomes singular.
The post-fault system trajectory hitting the impasse surface is regarded as one of the main causes for voltage collapse \cite{hiskens1989energy, venkatasubramanian1995local, praprost1996stability}.
Hence, the nature of an impasse surface is of importance to revealing the mechanism of voltage collapse.

\indent
Characterizing the impasse surface is a hard problem.
The existing results are mainly derived from simplified systems and only focus on the role of static loads.
For instance, the power system studied in \cite{oluic2018nature} is assumed to have one ZIP load and all the other loads are of constant-impedance type.
Those constant-impedance loads have no contribution to the impasse surface and are absorbed in the power network as shunt components, and the impact of the parameters of the single ZIP load on impasse surface is elaborated.
Hiskens and Hill \cite{hiskens1989energy} studied a four-bus system containing a single static load and proved that this specific system can avoid the impasse surface if the active power load is of constant-impedance type and the exponent of reactive power load is not less than one.
This also relates to \cite{lesieutre1999existence} that confirms the solvability of power system algebraic equation when active and reactive power load exponents are all greater than one.
In \cite{hiskens1990phd}, Hiskens extended the condition in \cite{hiskens1989energy} to a system containing two interconnected static loads and further conjectured that the result should also be valid for generic systems.
If true, this conjecture will provide an important class of load parameters that avoids the impasse surface; however, it remains unproved for decades.

\indent
To deepen the understanding of voltage collapse, the analysis of impasse surface needs to be extended to power systems with both static and dynamical loads. In addition, a general power network structure should be considered as the network structure is also crucial to system dynamics \cite{song2017cutset}, which fails to be captured by simplified system models.
A major obstacle in extending the existing methods is that they adopt certain assumptions or simplifications to obtain low-dimensional problem descriptions (e.g., scalar quadratic equations \cite{oluic2018nature, hiskens1989energy} or equations of up to 4$\times$4 matrices \cite{hiskens1990phd}) and derive explicit expressions for the spectrum or determinant of the algebraic Jacobian.
However, these tools are not applicable to generic cases with high-dimensional matrices, where the explicit solutions are unavailable.

\indent
In this paper, we develop an admittance matrix-based characterization for the impasse surface of DAE models of generic power systems by more advanced matrix analysis specific to the features of the algebraic Jacobian.
The following three aspects of contributions are made.

\indent
1) A necessary condition for a system trajectory hitting the impasse surface is established (see Theorem \ref{thmimpassenetw}).
It applies to a generic power system with multiple synchronous generators, static loads, induction motor loads, and a lossy power network.
This condition is in terms of admittance matrices regarding the effects of power network, generators and loads.
It carries clear network structural information and elaborates how the interactions between generators, loads and power network induce or prevent voltage collapse.
It also motivates an early indicator of voltage collapse to trigger corrective control.

\indent
2) Based on Theorem \ref{thmimpassenetw}, we manage to identify a class of static load parameters that make the system avoid the impasse surface (see Theorem \ref{thmimpasseexponent}).
This result proves the conjecture in \cite{hiskens1990phd} and has an even wider applicability.

\indent
3) We further show by eigen-analysis that inductive compensation to the power network has a positive effect on preventing voltage collapse, while capacitive compensation does the opposite, which is confirmed by simulation.

\indent
The remainder of the paper is organized as follows.
The DAE model of power systems is formulated in Section \ref{secformu}.
A new characterization of the impasse surface is given in Section \ref{secimpasse}. The obtained results are illustrated by simulation in Section \ref{seccase}.
Section \ref{secconclu} makes a conclusion and future prospect.

\indent
\textit{Notations:} The set of real numbers and complex numbers are denoted by $\mbR$ and $\mbC$, respectively.
The notation $\bm{x}=[x_i]\in \mbC^{p}$ denotes a vector, $\bm{x}=\diag{x_i}\in \mbC^{p\times p}$ denotes a diagonal matrix, and $\bm{I}_p\in \mbR^{p\times p}$ denotes an identity matrix.
In variation of the usual notation, the italic $j$ denotes a numbering index, while the upright j denotes the square root of -1.

\section{Power system differential-algebraic model}\label{secformu}
Consider a power system with $n$ buses coupled via a connected power network.
Suppose $g\leq n$ of the buses connect synchronous generators, and these buses are called the generator terminal buses.
Note that a generator can be modeled as an internal voltage source linking the corresponding terminal bus in the power network via a subtransient impedance (see Section~\ref{secformugen} for the details).
The power network is then augmented with $g$ buses and $g$ lines that represent the generator internal voltages and subtransient impedances.
Thus, there are totally $n+g$ buses with the addition of these ``virtual'' buses.
Let $\mV_G$ be the set of generator internal buses, $\mV_L$ be the set of remaining buses, $\mV_t\subseteq\mV_L$ be the set of buses connecting generator terminals and possibly loads, and $\mV_L\backslash\mV_t$ be the set of buses connecting loads only.
In the following we will formulate the system dynamical model by combining the models of generators, loads and network.

\subsection{Synchronous generator model}\label{secformugen}
We adopt a general high-order model for synchronous generators that includes the subtransient dynamics along $d$-axis and $q$-axis and possibly an excitation system (e.g., automatic voltage regulator).
There are several representative generator models considering subtransient dynamics, such as the Sauer-Pai's model, Marconato's model and Anderson-Fouad's model~\cite{milano2010power}.
We do not concern much about the details of generator differential equations since the impasse surface is only concerned with the algebraic equations and algebraic variables.
In the following, we will present the generator model in a compact form for simplicity.

\indent
Let $\theta_i, V_i$ denote the phase angle and voltage magnitude of bus $i\in \mV_L$.
For the generator connecting to terminal bus $i\in \mV_t\subseteq\mV_L$, its state variables consist of the rotor angle $\delta_i$, rotor speed $\omega_i$, transient $d$-axis and $q$-axis voltages $E_{di}^{\prime}, E_{qi}^{\prime}$, subtransient $d$-axis and $q$-axis voltages $E_{di}^{\prime\prime}, E_{qi}^{\prime\prime}$, and possibly excitation system variable $\bm{x}_{fi}$.
The algebraic variables associated with this generator are $\theta_i, V_i$.
In general, the generator dynamics can be described by
\begin{equation}\label{gensixth}
\begin{split}
   \dot{\bm{x}}_{gi} = \bm{f}_{gi}(\bm{x}_{gi}, \theta_i, V_i)
\end{split}
\end{equation}
where $\bm{x}_{gi}=\begin{bmatrix}
                    \delta_i & \omega_i & E_{di}^{\prime} & E_{qi}^{\prime} & E_{di}^{\prime\prime} & E_{qi}^{\prime\prime} &  \bm{x}_{fi}^T \\
                  \end{bmatrix}^T$
collects the state variables.
Note that \eqref{gensixth} is also dependent on some important parameters including the armature resistance $r_{ai}$, synchronous reactances $x_{di}, x_{qi}$, transient reactances $x_{di}^{\prime}, x_{qi}^{\prime}$ and subtransient reactances $x_{di}^{\prime\prime}, x_{qi}^{\prime\prime}$.
But these parameters are constant coefficients in the model so that we do not explicitly express them in \eqref{gensixth}.

\indent
In addition, we adopt the assumption below for the synchronous generators.

\begin{assumption}\label{assump1}
   For each generator, the $d$-axis subtransient reactance $x_{di}^{\prime\prime}$ is equal to the $q$-axis subtransient reactance $x_{qi}^{\prime\prime}$.
\end{assumption}

\begin{figure}[!h]
  \centering
  \includegraphics[width=3.3in]{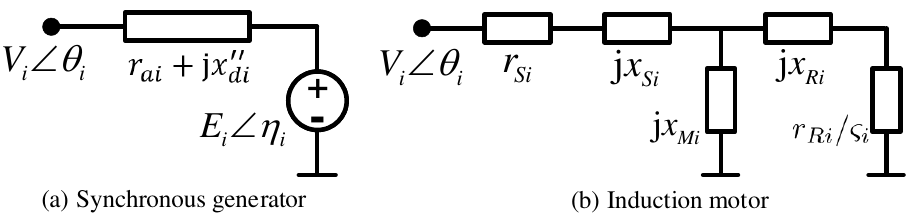}
  \caption{The equivalent circuit for a generator and an induction motor.}
  \label{figgenmot}
\end{figure}

\indent
The difference between $x_{di}^{\prime\prime}$ and $x_{qi}^{\prime\prime}$ is called subtransient saliency.
For a generator with a damper winding in both the $d$-axis and the $q$-axis, the screening effect in both axes is similar and subtransient saliency is negligible \cite{machowski2008power}.
So Assumption~\ref{assump1} is reasonable and has been commonly used to simplify generator modeling~\cite{kundur1994power}.
With this assumption, we can obtain a neat equivalent circuit for the generator as shown in Fig. \ref{figgenmot}(a) with two features.
First, the impedance linking the generator internal bus and terminal bus is the subtransient impedance $r_{ai}+\upj x_{di}^{\prime\prime}$.
Second, the phase angle and voltage magnitude of the generator internal bus, say $\eta_i$ and $E_i$, are functions of the state variables only, i.e., $\eta_i= \eta_i(\bm{x}_{gi})$ and $E_i= E_i(\bm{x}_{gi})$.
The expressions of these two functions vary with the generator model (e.g., see \cite{sauer1998power, milano2010power}).
For instance, we have the following equation when adopting the Sauer-Pai's model \cite{milano2010power}
\begin{equation}\label{sauerpaimodel}
\begin{split}
   0 &= r_{ai}i_{qi} + x_{di}^{\prime\prime}i_{di} + V_{qi} - \gamma_{di}E_{qi}^{\prime} - (1-\gamma_{di})E_{qi}^{\prime\prime}\\
   0 &= r_{ai}i_{di} - x_{qi}^{\prime\prime}i_{qi} + V_{di} - \gamma_{qi}E_{di}^{\prime} - (1-\gamma_{qi})E_{di}^{\prime\prime}
\end{split}
\end{equation}
where $V_i\angle{\theta_i}=V_{di}+\upj V_{qi}$, $i_{di}+\upj i_{qi}$ denotes the generator current and $\gamma_{di}, \gamma_{qi}$ are determined by $x_{di}^{\prime}, x_{qi}^{\prime}, x_{di}^{\prime\prime}, x_{qi}^{\prime\prime}$.
If $x_{di}^{\prime\prime}=x_{qi}^{\prime\prime}$, then we can set
$E_i\angle{\eta_i}= \gamma_{qi}E_{di}^{\prime}+(1-\gamma_{qi})E_{di}^{\prime\prime} + \upj \gamma_{di}E_{qi}^{\prime} +\upj (1-\gamma_{di})E_{qi}^{\prime\prime}$
as a function of state variables only and obtain from \eqref{sauerpaimodel} that $E_i\angle{\eta_i}=V_i\angle{\theta_i}+(r_{ai}+\upj x_{di}^{\prime\prime})(i_{di}+\upj i_{qi})$, which is consistent with Fig.~\ref{figgenmot}(a).
Since $\eta_i$ and $E_i$ are independent of algebraic variables, they can be treated as constants in the differentiation with respect to algebraic variables, which will bring convenience to the analysis of impasse surface in Section \ref{secimpasse}.

\subsection{Load model}
We adopt the composite load model consisting of static components and induction motor components, which is a common model for stability analysis concerning voltage dynamics \cite{milanovic2013international, zhang2015ambient}.
The static load component at bus $i\in \mV_L$ is described by the exponential terms
\begin{equation}\label{staticload}
\begin{split}
      P_{si}(V_i) &= P_{si}^0 V_i^{\alpha_i}\\
      Q_{si}(V_i) &= Q_{si}^0 V_i^{\beta_i}
\end{split}
\end{equation}
where $\alpha_i, \beta_i$ denote the active and reactive power load exponent, and $P_{si}^0, Q_{si}^0$ denote the rated active and reactive power load.

\indent
Next, we model the induction motors.
Different from the synchronous generators whose rotors run exactly at the system frequency at an equilibrium point, the steady-state rotor frequency of an induction motor is lower than the system frequency and its power consumption depends on the motor slip.
In this paper, we adopt the third-order model \cite{milano2010power} to capture both electromechanical and electromagnetic dynamics of induction motors.
The equivalent circuit for the third-order induction motor model is depicted in Fig. \ref{figgenmot}(b), where $\varsigma_i$ denotes the slip, $r_{Si}, x_{Si}$ denotes the motor stator resistance and reactance, $r_{Ri}, x_{Ri}$ denotes the cage rotor resistance and reactance, and $x_{Mi}$ denotes the magnetization reactance.
Then, the equivalent admittance of the motor circuit is
\begin{equation}\label{Ymi}
\begin{split}
      Y_{mi}^{\upeq}(\varsigma_i) = \Big( r_{Si}+\upj x_{Si} + \frac{\upj x_{Mi}(r_{Ri}/\varsigma_i+\upj x_{Ri})}{r_{Ri}/\varsigma_i+\upj (x_{Ri}+x_{Mi})}\Big)^{-1}
\end{split}
\end{equation}
and hence the motor load consumption at bus $i$ is
\begin{equation}\label{motorload}
\begin{split}
      P_{mi}(\varsigma_i, V_i) &= G_{mi}^{\upeq}(\varsigma_i)V_i^2\\
      Q_{mi}(\varsigma_i, V_i) &= -B_{mi}^{\upeq}(\varsigma_i)V_i^2
\end{split}
\end{equation}
where $G_{mi}^{\upeq}$ and $B_{mi}^{\upeq}$ respectively denote the real part and imaginary part of $Y_{mi}^{\upeq}$.

\indent
In addition, we have the following differential equations that describe the motion and internal voltage of an induction motor
\begin{equation}\label{motorDE}
\begin{split}
      \dot{\bm{x}}_{mi} = \bm{f}_{mi}(\bm{x}_{mi}, \theta_i, V_i)
\end{split}
\end{equation}
where $\bm{x}_{mi}=\begin{bmatrix}
                    \varsigma_i & e_{di}^{\prime} & e_{qi}^{\prime} \\
                  \end{bmatrix}^T$ with $e_{di}^{\prime}, e_{qi}^{\prime}$ being the $d$-axis and $q$-axis voltage behind the the stator resistance $r_{Si}$.
Since we will focus on the algebraic equation, again we omit the explicit form of \eqref{motorDE} and refer to \cite{milano2010power} for the details.

\subsection{Power network model \& power flow equation}
Without loss of generality, $\mV_t$, $\mV_L$ and $\mV_G$ are numbered as $\mV_t=\{1,...,g\}$, $\mV_L=\{1,...,n\}$ and $\mV_G=\{n+1,...,n+g\}$.
Let $\bm{Y}_{\bus}=[Y_{ij}]\in\mbC^{n\times n}$ be the power network admittance matrix among $\mV_L$, which is defined by
\begin{equation}
\begin{split}
     Y_{ij} =
       \left\{
        \begin{array}{ll}
            y_{i0} + \sum_{j=1,j\neq i}^{n} y_{ij} ,~i=j\in\mV_L \\
            -y_{ij},~i\neq j,~i,j\in\mV_L
        \end{array}
        \right.
\end{split}
\end{equation}
where $y_{ij}\in\mbC$ denotes the admittance of line $(i,j)$, $y_{ij}=0$ if bus $i$ and bus $j$ are not directly connected;
$y_{i0}\in\mbC$ denotes the shunt component at bus $i$ such as the line charging capacitance.
The matrix $\bm{Y}_{\bus}$ is commonly used in power system steady-state analysis where the generator internals $\mV_G$ are not considered~\cite{kundur1994power}.
For studying the DAE model, we also need to introduce the augmented admittance matrix including $\mV_G$, say $\widetilde{\bm{Y}}=[\widetilde{Y}_{ij}]\in\mbC^{(n+g)\times(n+g)}$, which takes the form
\begin{equation}\label{Yaug2Y}
\begin{split}
   \widetilde{\bm{Y}} =
   \begin{bmatrix}
     \bm{Y}_{\bus}+\bm{Y}_{\gen} & \bm{Y}_{LG} \\
     \bm{Y}_{LG}^T & \bm{Y}_{GG} \\
   \end{bmatrix}
\end{split}
\end{equation}
where $\bm{Y}_{GG}=\diag{y_i^{gs}}\in\mbC^{g\times g}$ with $y_i^{gs}=(r_{si}+\upj x_{di}^{\prime\prime})^{-1}$, $\forall i\in\mV_t$;
$\bm{Y}_{\gen}=\diag{Y_{gi}}\in\mbC^{n\times n}$ with $Y_{gi}=y_i^{gs}$ if $i\in\mV_t$ and $Y_{gi}=0$ if $i\in\mV_L\backslash\mV_t$;
$\bm{Y}_{LG}\in\mbC^{n\times g}$ is defined such that $[\bm{Y}_{LG}]_{ij}=-y_i^{gs}$ if bus $i\in\mV_t$ and bus $n+j$ is the corresponding generator internal bus, and $[\bm{Y}_{LG}]_{ij}=0$ otherwise.
We will use $\widetilde{G}_{ij}, G_{ij}, G_{gi}$ (or $\widetilde{B}_{ij}, B_{ij}, B_{gi}$) denote the real (or imaginary) parts of $\widetilde{Y}_{ij}, Y_{ij},Y_{gi}$, respectively.

\indent
To obtain a neat expression for the power flow equation, we henceforth use $\theta_i, V_i$ to denote the phase angle and voltage magnitude for any bus $i\in\mV_L\cup\mV_G$.
For notation consistency, for each generator internal bus $j\in\mV_G$ and its associated terminal bus $i\in\mV_t$, we set
\begin{equation}\label{EdEq2V}
\begin{split}
   V_j = E_i(\bm{x}_{gi}),~\theta_j = \eta_i(\bm{x}_{gi}).
\end{split}
\end{equation}
Then, the power balance at each bus $i\in\mV_L$ can be described by the power flow equation\footnote{The generator internal buses $\mV_G$ only achieve power balance at an equilibrium. In that case, \eqref{gensixth} degenerates to an algebraic equation regarding power balance. So there is no specific power flow equation for bus $i\in\mV_G$.}
\begin{equation}\label{pfequ}
\begin{split}
   0 = g_{pi} =&  V_i^2\widetilde{G}_{ii} + \sum\nolimits_{j=1,j\neq i}^{n+g} V_iV_j|\widetilde{Y}_{ij}|\sin(\theta_{ij}-\varphi_{ij})\\ &
   + P_{si}(V_i) + P_{mi}(\varsigma_i, V_i) \\
   0 = g_{qi} =& -V_i^2\widetilde{B}_{ii} - \sum\nolimits_{j=1,j\neq i}^{n+g} V_iV_j|\widetilde{Y}_{ij}|\cos(\theta_{ij}-\varphi_{ij})\\ &
   + Q_{si}(V_i) + Q_{mi}(\varsigma_i, V_i)
\end{split}
\end{equation}
where $\theta_{ij}$ is defined as $\theta_{ij}=\theta_i-\theta_j$;
and $\varphi_{ij}=-\tan^{-1}(\widetilde{G}_{ij}/\widetilde{B}_{ij})$ is the phase shift caused by line loss (we set $\varphi_{ij}=0$ if $\widetilde{Y}_{ij}=0$).
In fact, the power network does have electromagnetic dynamics; however, these dynamics decay much faster than generator electromechanical swings and load behaviors.
Hence, it is reasonable to describe the power network by \eqref{pfequ}, which is common in power system communities \cite{van2007voltage}.
Moreover, the right hand side of \eqref{pfequ}, which is the core of following analysis, is indeed the coupling terms in general Kuramoto oscillators with heterogeneous edge weights and phase shifts \cite{dorfler2012synchronization}. The obtained results will also be useful to the study of Kuramoto oscillator dynamics.

\indent
To sum up, let $\bm{x}, \bm{y}$ be the vectors of state variables and algebraic variables, respectively, where $\bm{x}$ collects $\bm{x}_{gi}$, $\forall i\in\mV_t$ and $\bm{x}_{mi}$, $\forall i\in\mV_L$, and $\bm{y}$ collects $\bm{\theta}=[\theta_i]\in\mbR^n, \bm{V}=[V_i]\in\mbR^n$, $\forall i\in\mV_L$.
Then, power system dynamics can be described by the following DAE model in the compact form
\begin{subequations}\label{DAE}
\begin{align}
   \dot{\bm{x}} &= \bm{f}(\bm{x}, \bm{y}) \label{DAEdiff}  \\
   \bm{0} &= \bm{g}(\bm{x}, \bm{y})   \label{DAEalg}
\end{align}
\end{subequations}
where the differential equations \eqref{DAEdiff} consist of \eqref{gensixth} and \eqref{motorDE}, and the algebraic equations \eqref{DAEalg} consist of \eqref{pfequ}.
Note that the algebraic variables do not include $\theta_j, V_j$, $\forall j\in\mV_G$ as they can be substituted by the associated state variables using~\eqref{EdEq2V}.

\section{Characterizing impasse surface}\label{secimpasse}
\subsection{Impasse surface \& voltage collapse}
Let us first introduce a matrix closely linked to the impasse surface, namely the algebraic Jacobian $\Jalg(\bm{x}, \bm{y})=\partial \bm{g}(\bm{x}, \bm{y})/\partial \bm{y}$, which by \eqref{pfequ} can be further expanded by
\begin{equation}\label{pfJacobian}
\begin{split}
   \Jalg(\bm{x}, \bm{y})=
   \begin{bmatrix}
     \frac{\partial \bm{g}_p}{\partial \bm{\theta}} & \frac{\partial \bm{g}_p}{\partial \bm{V}} \\
     \frac{\partial \bm{g}_q}{\partial \bm{\theta}} & \frac{\partial \bm{g}_q}{\partial \bm{V}} \\
   \end{bmatrix}\in \mbR^{2n\times 2n}.
\end{split}
\end{equation}
The entries of $\frac{\partial \bm{g}_p}{\partial \bm{\theta}}, \frac{\partial \bm{g}_p}{\partial \bm{V}}, \frac{\partial \bm{g}_q}{\partial \bm{\theta}}, \frac{\partial \bm{g}_q}{\partial \bm{V}}\in \mbR^{n\times n}$, respectively denoted by $\frac{\partial g_{pi}}{\partial \theta_j}, \frac{\partial g_{pi}}{\partial V_j}, \frac{\partial g_{qi}}{\partial \theta_j}, \frac{\partial g_{qi}}{\partial V_j}$, $i,j=1,...,n$, take values as
\begin{equation}\label{Jpt}
\begin{split}
     \frac{\partial g_{pi}}{\partial \theta_j} &=
       \left\{
        \begin{array}{ll}
           \sum_{j=1,j\neq i}^{n+g} V_iV_j|\widetilde{Y}_{ij}|\cos(\theta_{ij}-\varphi_{ij}),~i=j \\
           -V_iV_j|\widetilde{Y}_{ij}|\cos(\theta_{ij}-\varphi_{ij}),~i\neq j.
        \end{array}
        \right.\\
     \frac{\partial g_{pi}}{\partial V_j} &=
     \left\{
        \begin{array}{ll}
           \sum_{j=1,j\neq i}^{n+g} V_j|\widetilde{Y}_{ij}|\sin(\theta_{ij}-\varphi_{ij}) \\
           ~~~+ 2V_i(\widetilde{G}_{ii}+G_{mi}^{\upeq}) + \alpha_i P_{si}^0 V_i^{\alpha_i-1},~i=j \\
           V_i|\widetilde{Y}_{ij}|\sin(\theta_{ij}-\varphi_{ij}),~i\neq j.
        \end{array}
        \right.\\
     \frac{\partial g_{qi}}{\partial \theta_j} &=
      \left\{
        \begin{array}{ll}
           \sum_{j=1,j\neq i}^{n+g} V_iV_j|\widetilde{Y}_{ij}|\sin(\theta_{ij}-\varphi_{ij}),~i=j \\
           -V_iV_j|\widetilde{Y}_{ij}|\sin(\theta_{ij}-\varphi_{ij}),~i\neq j.
        \end{array}
        \right.\\
     \frac{\partial g_{qi}}{\partial V_j} &=
       \left\{
        \begin{array}{ll}
           -\sum_{j=1,j\neq i}^{n+g} V_j|\widetilde{Y}_{ij}|\cos(\theta_{ij}-\varphi_{ij}) \\
           ~~ -2V_i(\widetilde{B}_{ii}+B_{mi}^{\upeq}) + \beta_i Q_{si}^0 V_i^{\beta_i-1},~i=j \\
           -V_i|\widetilde{Y}_{ij}|\cos(\theta_{ij}-\varphi_{ij}),~i\neq j.
        \end{array}
        \right.
\end{split}
\end{equation}
When we are differentiating the algebraic equations with respect to algebraic variables to obtain $\Jalg$, the terms $\theta_j, V_j$, $j\in\mV_G$ are regarded as constants as they are functions of the state variables by~\eqref{EdEq2V} and independent of the algebraic variables.
Nevertheless, the entries of $\Jalg$ depend on both the state variables and algebraic variables as the terms $\theta_j, V_j$, $j\in\mV_G$, which are functions of the state variables, still appear in \eqref{Jpt}.

\indent
Then, the impasse surface is defined in terms of $\Jalg$ below.

\begin{definition}[\cite{hiskens1989energy}]
   The impasse surface of system \eqref{DAE} consists of the set of points $\textup{IS}=\{(\bm{x}, \bm{y})|~\textup{det}\{\Jalg(\bm{x}, \bm{y})\}=0\}$.
\end{definition}

\indent
As mentioned before, the DAE model \eqref{DAE} ignores those fast dynamics of a power system such as the network electromagnetic dynamics.
Nevertheless, the impasse surface of the DAE model can be used to interpret some physical behaviours for the power system.
Let $(\bm{x}(t), \bm{y}(t))$ be a trajectory of \eqref{DAE}.
When the trajectory hits the impasse surface, the DAE model losses causality since the algebraic variable $\bm{y}$ can no longer be predicted by state variable $\bm{x}$ from the relation $\frac{\partial \bm{g}}{\partial \bm{x}}\Delta\bm{x}+\Jalg\Delta\bm{y}=\bm{0}$ with a singular $\Jalg$. The time-domain simulation of the DAE model fails to continue afterwards.
On the other hand, a more detailed model of the power system is in the form of pure differential equations (DEs), which is obtained by replacing the algebraic equations \eqref{DAEalg} by DEs with very small time constants to capture those ignored fast dynamics \cite{sauer1987reduced}.
This DE model never fails to continue at any point, but it will have undesirable behaviors when its associated DAE model hits the impasse surface.
It is widely observed that the states (e.g., bus voltages) of the DE model have a rapid movement when its associated DAE model hits the impasse surface \cite{sastry1981jump, chua1989impasse1}, which corresponds to short-term voltage collapse.

\indent
Moreover, the analytical study of impasse surfaces commonly takes the following assumption.

\begin{assumption}\label{assump2}
   The voltage magnitudes of all buses are non-zero along the system trajectory $(\bm{x}(t), \bm{y}(t))$.
\end{assumption}

\indent
By \eqref{Jpt}, the algebraic Jacobian has its $i$-th column being zero and hence is singular if the voltage of bus $i$ becomes zero along the trajectory.
Note that zero voltages occur only when the system undergoes a purely metallic short-circuit fault (i.e., the fault impedance is strictly zero), which is rare in practice. Also the zero voltages in this case already give a clear indication of collapse.
So we focus on the nontrivial case that voltage collapse occurs when bus voltages are still away from zero, which is harder to detect and of more interest.

\indent
In case that all loads are purely static (i.e., $Y_{mi}^{\upeq}=0$, $\forall i\in\mV_L$), Hiskens proposed the following conjecture for an impasse surface in \cite{hiskens1990phd}.

\begin{conjecture}[\cite{hiskens1990phd}]\label{conjimpasse}
     Suppose Assumption \ref{assump1} and Assumption \ref{assump2} hold.
     Consider a DAE system \eqref{DAE} with all loads being purely static.
     The system trajectories never encounter the impasse surface if the following conditions are all satisfied:
     \begin{enumerate}
       \item Generator circuit: $G_{gi}=0$, $B_{gi}<0$, $\forall i\in\mV_t$;
       \item Power network: $G_{ij}=0$, $B_{ij}>0$, $B_{ii}= -\sum\nolimits_{j=1,j\neq i}^n B_{ij}$, $\forall i,j\in\mV_L$;
       \item Active power load: $P_{si}^0\geq 0$, $\alpha_i=2$, $\forall i\in\mV_L$;
       \item Reactive power load: $Q_{si}^0\geq0$, $\beta_i\geq 1$, $\forall i\in\mV_L$.
     \end{enumerate}
\end{conjecture}

\indent
This conjecture can provide a class of load parameters that prevent system trajectories from hitting the impasse surface once it is confirmed true.
In the following, we will link the impasse surface to admittance matrices with new insights into voltage collapse.
Further, our analysis proves Conjecture~\ref{conjimpasse}.

\subsection{Theoretical results and physical interpretations}\label{sectheory}
We begin the analysis by defining the equivalent conductances and susceptances of static loads as follows.

\begin{definition}\label{defYsi}
   For each bus $i\in\mV_L$, define
   \begin{equation}\label{Ysi}
   \begin{split}
     G_{si}^{\upeq}(t) &= P_{si}(t)/V_i^2(t)\\
     B_{si}^{\upeq}(t) &= -Q_{si}(t)/V_i^2(t)
   \end{split}
   \end{equation}
   as the equivalent conductance and equivalent susceptance of the static load at bus $i$ at time $t$, respectively.
\end{definition}

\indent
This definition has straightforward physical meanings.
At any time $t$, if we replace the static load at bus $i$ by the shunt admittance $G_{si}^{\upeq}(t)+\upj B_{si}^{\upeq}(t)$, then its power consumption is exactly $P_{si}(t)+\upj Q_{si}(t)$.
In the following, we will substitute the equivalent conductances and susceptances of induction motors and static loads $G_{mi}^{\upeq},B_{mi}^{\upeq}, G_{si}^{\upeq},B_{si}^{\upeq}$ into the algebraic Jacobian entries and derive a novel condition on the impasse surface.

\indent
For the convenience of presenting our results, we introduce the following admittance matrices regarding the induction motors and static loads $\bm{Y}_{\mot}^{\upeq}(t)=\diag{Y_{mi}^{\upeq}(\varsigma_i(t))}, \bm{G}_{\stat}^{\upeq}(t)=\diag{G_{si}^{\upeq}(t)}, \bm{B}_{\stat}^{\upeq}(t)=\diag{B_{si}^{\upeq}(t)},\bm{\alpha}=\diag{\alpha_i},\bm{\beta}=\diag{\beta_i}\in\mbR^{n\times n}$.
Then, we are ready to state the following theorem (the proof is given in Appendix).

\begin{theorem}\label{thmimpassenetw}
     Suppose Assumption \ref{assump1} and Assumption \ref{assump2} hold.
     The trajectory of DAE system \eqref{DAE} encounters the impasse surface at time $t$ only if
     \begin{equation}\label{impassenetwineq}
     \begin{split}
        \sigma_{\min}(\bm{Y}_1(t))&\leq \max_{i\in\mV_L}~\big|(1-\frac{\alpha_i}{2})G_{si}^{\upeq}(t) + \upj (1-\frac{\beta_i}{2})B_{si}^{\upeq}(t)\big|
     \end{split}
     \end{equation}
     where $\sigma_{\min}$ denotes the minimum singular value and
     \begin{equation}\label{Y1}
     \begin{split}
        \bm{Y}_1(t) &= \bm{Y}_{\bus}+\bm{Y}_{\gen}+\bm{Y}_{\mot}^{\upeq}+\frac{1}{2}\bm{\alpha}\bm{G}_{\stat}^{\upeq}+\upj \frac{1}{2}\bm{\beta}\bm{B}_{\stat}^{\upeq}.
     \end{split}
     \end{equation}
\end{theorem}

\indent
Theorem \ref{thmimpassenetw} establishes a necessary condition for a system trajectory hitting the impasse surface by the admittance matrices of the power network, generator equivalent circuits, static loads and induction motors.
This result has wide applicability as it allows a generic modeling for generators, loads and power network, which sheds new light on the role of these system components in inducing voltage collapse.
It generalizes the results in \cite{oluic2018nature} which focuses on the parameters of a single load with the other loads being constant impedances.

\indent
Observing the admittance terms in \eqref{impassenetwineq}, $\bm{Y}_{\bus}$ refers to the coupling among buses $\mV_L$, $\bm{Y}_{\gen}$ refers to the coupling between generator internals and terminals, $\bm{Y}_{\mot}^{\upeq}$ can be regarded as the coupling between the power network and induction motors, and $\bm{G}_{\stat}^{\upeq}, \bm{B}_{\stat}^{\upeq}$ represents the effect of static loads.
When the static loads are of constant power type (i.e., $\bm{\alpha}=\bm{\beta}=\bm{0}$), the terms with respect to static loads vanish in the left-hand-side of \eqref{impassenetwineq}.
In this case, Theorem \ref{thmimpassenetw} leads to an intuitive interpretation of short-term voltage stability, i.e., the system trajectory avoids hitting the impasse surface if the coupling between power network, generators and motors is sufficiently strong to ``prevail over'' the effect of static loads.
Further, in generic cases with non-zero load exponents, the effect of static loads contributes to both sides of inequality \eqref{impassenetwineq}.

\indent
Theorem \ref{thmimpassenetw} also coincides with the intuition that a low voltage level must occur during collapse.
At a ``healthy'' state where $V_i\simeq$1.0 p.u., $\bm{Y}_{\bus}$ and $\bm{Y}_{\gen}$ are much greater than the other terms relating to equivalent load admittances so that inequality \eqref{impassenetwineq} is not satisfied.
On the other hand, a severe voltage decline caused by a disturbance (e.g., short-circuit fault) significantly increases the equivalent load admittances, which makes it possible to satisfy \eqref{impassenetwineq} and eventually induces voltage collapse.

\indent
Further, in visualizing the process of voltage collapse, it is convenient to define the index $\Ivs(t)$ as the ratio of the left-hand-side to right-hand-side of \eqref{impassenetwineq}.
According to Theorem~\ref{thmimpassenetw}, $\Ivs(t)$ being less than one is a necessary condition for hitting the impasse surface, which means the actual time when system trajectory hits the impasse surface must be later than the time when $\Ivs(t)$ is below one.
Hence, the index $\Ivs(t)$ is a dynamic indicator that can provide an early warning of voltage collapse for triggering corrective control, an example of which will be shown in the case study.
By comparison, a necessary and sufficient condition for hitting the impasse surface will not leave any time for corrective control. It shows that Theorem~\ref{thmimpassenetw}, which inevitably has conservativeness as a necessary condition, does have some merits in control application.

\indent
Next, we move to the special case where all loads are purely static (i.e., $Y_{mi}^{\upeq}(t)=0$, $\forall i\in\mV_L$), which is commonly studied in the literature \cite{hiskens1989energy, lesieutre1999existence, oluic2018nature}.
In this case we have the following theorem regarding the impact of load exponents (the proof is given in Appendix).

\begin{theorem}\label{thmimpasseexponent}
     Suppose Assumption \ref{assump1} and Assumption \ref{assump2} hold.
     Consider a DAE system \eqref{DAE} with all loads being purely static.
     The system trajectories never encounter the impasse surface if the following conditions are all satisfied:
     \begin{enumerate}
       \item Generator circuit: $G_{gi}\geq 0$, $B_{gi}<0$, $\forall i\in\mV_t$;
       \item Power network: $G_{ij}=0$, $B_{ij}>0$, $B_{ii}= -\sum\nolimits_{j=1,j\neq i}^n B_{ij}$, $\forall i,j\in\mV_L$;
       \item Active power load: $P_{si}^0=0$ or $\alpha_i=2$, $\forall i\in\mV_L$;
       \item Reactive power load: $Q_{si}^0\geq0$, $\beta_i\geq 1$, $\forall i\in\mV_L$.
     \end{enumerate}
\end{theorem}

\indent
We further interpret the conditions in Theorem \ref{thmimpasseexponent}.
The generator circuit condition is general and trivial.
The power network condition is a reasonable approximation for the situation in high-voltage transmission systems, where the lines are inductive (i.e., $B_{ij}>0$) with negligible conductance (i.e., $G_{ij}=0$) and the charging capacitance is negligible compared to line susceptance (i.e., $B_{ii}= -\sum\nolimits_{j=1,j\neq i}^n B_{ij}$).
The load condition approximates such an operating scenario that the active power loads are very small or behave like constant impedances.
It indicates that active power loads have no contribution to the right-hand-side of \eqref{impassenetwineq}, or in other words, the reactive power loads are the dominant factor.
This also corresponds to a typical scenario for the voltage stability issues in transmission systems where the voltages are more strongly coupled with reactive powers than active powers.

\indent
In general, Theorem~\ref{thmimpasseexponent} implies that the system is much less likely to suffer voltage collapse when the active power loads are constant impedances and reactive power load exponents are no less than one.
Particularly, it proves Conjecture~\ref{conjimpasse} with even more relaxed conditions.
Unlike Conjecture~\ref{conjimpasse}, Theorem~\ref{thmimpasseexponent} still holds if:
1) $G_{gi}> 0$ which corresponds to a lossy generator circuit; and
2) $P_{si}^0<0$ and $\alpha_i=2$ which corresponds to a ``negative load'' case that could be the result of demand-side management or high penetration of renewable energy.

\indent
Theorem~\ref{thmimpasseexponent} also relates to some existing findings on the impact of load exponents.
For instance, it is observed in \cite{oluic2018nature} that it is highly difficult to find an event of hitting the impasse surface in case of constant-current loads (i.e., $\alpha_i=\beta_i=1$), which only occurs at an unrealistically heavy load level.
It is proved in \cite{lesieutre1999existence} that the algebraic equation \eqref{DAEalg} exhibits at least one solution if $\alpha_i>1$ and $\beta_i>1$, $\forall i\in\mV_L$.
The solution existence of algebraic equation almost indicates the non-singularity of algebraic Jacobian, except when the algebraic equation has a unique solution in some critical situations.
Theorem~\ref{thmimpasseexponent} is consistent with these results with new insights.

\subsection{Impact of shunt capacitor/inductor}\label{secshunt}
Shunt capacitors and shunt inductors are common devices for reactive power compensation and voltage regulation.
Based on the obtained theorems, this subsection carries out a qualitative analysis for the role of shunt capacitors and shunt inductors in voltage collapse.

\indent
The shunt devices can be regarded as a part of the power network. If bus $i$ installs a shunt capacitor/inductor, then it adds a term $\upj b_{i0}$ to the $(i,i)$-entry of $\bm{Y}_{\bus}$, where $b_{i0}>0$ implies capacitive compensation and $b_{i0}<0$ implies inductive compensation.
By Theorem \ref{thmimpassenetw}, the impact of $\upj b_{i0}$ on preventing/causing voltage collapse can be reflected by how it affects the matrix $\bm{Y}_1$ defined in \eqref{Y1}. If $\sigma_{\min}(\bm{Y}_1)$ is increased (or decreased) after adding $\upj b_{i0}$, then it implies that the system trajectory is less (or more) likely to hit the impasse surface, and hence a lower (or higher) risk of voltage collapse.

\indent
Before proceeding further, we make an approximation that the real part of $\bm{Y}_1$ is negligible compared to its imaginary part.
This can be justified by the usual case where the line conductances and equivalent load conductances are much smaller than line susceptances.
Then, we have $\bm{Y}_1=\upj \bm{B}_1$ where $\bm{B}_1= \bm{B}_{\bus}+\bm{B}_{\gen}+\bm{B}_{\mot}^{\upeq}+\frac{1}{2}\bm{\beta}\bm{B}_{\stat}^{\upeq}\in\mbR^{n\times n}$ with $\bm{B}_{\bus},\bm{B}_{\gen},\bm{B}_{\mot}^{\upeq}$ being the imaginary part of $\bm{Y}_{\bus},\bm{Y}_{\gen},\bm{Y}_{\mot}^{\upeq}$, respectively.
If $B_{ij}>0$ for $i\neq j$, $B_{ii}= -\sum\nolimits_{j=1,j\neq i}^n B_{ij}<0$, $B_{gi}<0$ and $Q_{si}^0\geq0$, which commonly holds in transmission systems, then $-\bm{B}_1$ is positive definite as it can be regarded as a graph Laplacian matrix with positive weighted lines and positive self-loops \cite{dorfler2013kron}.
It follows that $\sigma_{\min}(\bm{Y}_1)=\lambda_{\min}(-\bm{B}_1)$, where $\lambda_{\min}$ denotes the minimum eigenvalue.
Then, when adding a shunt inductor $\upj b_{i0}$ with $b_{i0}<0$ (or a capacitor $\upj b_{i0}$ with $b_{i0}>0$) to bus $i$, it decreases (or increases) $B_{ii}$ and hence increases (or decreases) the $i$-th main diagonal of $-\bm{B}_1$. Thus, by eigenvalue sensitivity analysis \cite{petersen2008matrix}, $\sigma_{\min}(\bm{Y}_1)=\lambda_{\min}(-\bm{B}_1)$ is increased (or decreased) after adding a shunt inductor (or capacitor) to bus~$i$.
This result indicates that inductive compensation has a positive effect on preventing voltage collapse while capacitive compensation does the opposite. An example will be given in the case study.

\section{Case study}\label{seccase}
\begin{figure}[!h]
  \centering
  \includegraphics[width=2.0in]{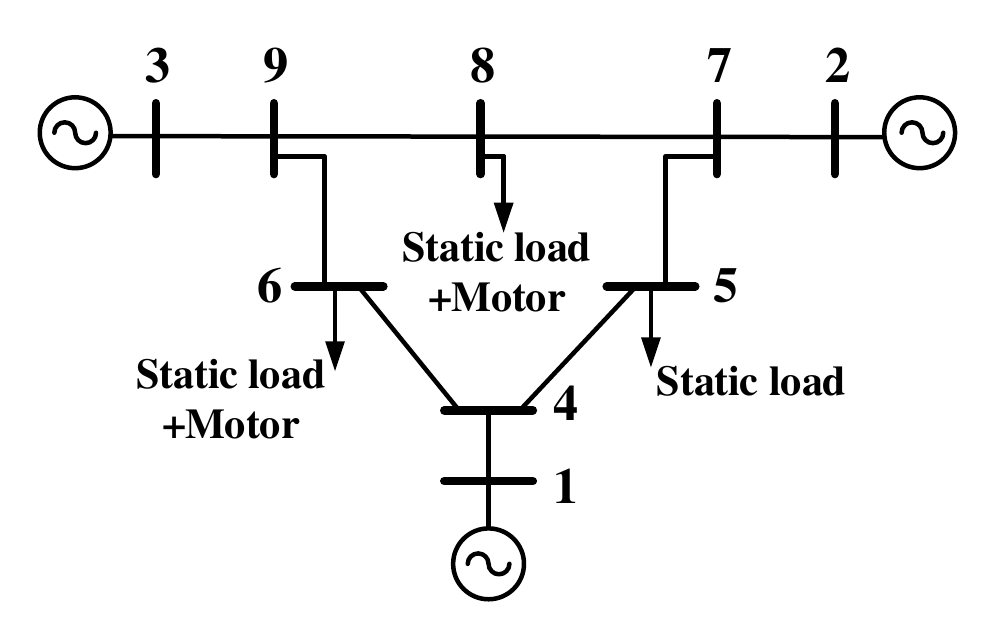}
  \caption{Diagram of the IEEE 9-bus system.}
  \label{fig9busdiagram}
\end{figure}

Take the IEEE 9-bus system to illustrate the obtained results.
The system diagram is given in Fig. \ref{fig9busdiagram}, where the generator internal buses are not displayed for simplicity.
In brief, bus 1, bus 2 and bus 3 are generator terminals and bus 5, bus 6 and bus 8 connect loads.
The load at bus 5 is purely static with $\alpha_5=0.1$, $\beta_5=0.6$. The loads at bus 6 and bus 8 consist of induction motors and static components with $\alpha_6=1.0$, $\beta_6=1.0$ and $\alpha_8=0.4$, $\beta_8=0.4$. The generators at bus 1 and bus 2 install the simplified IEEE Type DC1 excitor \cite{kundur1994power}, and the generator at bus 3 has no excitation system and keeps a constant field voltage.
We refer to \cite{case9data2019} for the detailed model built in PSAT \cite{milano2005open} format.

\begin{figure}[!h]
  \centering
  \subfigure[Voltage of bus 8.]{
  \label{figcase9voltage} %% label for first subfigure
  \includegraphics[width=1.67in]{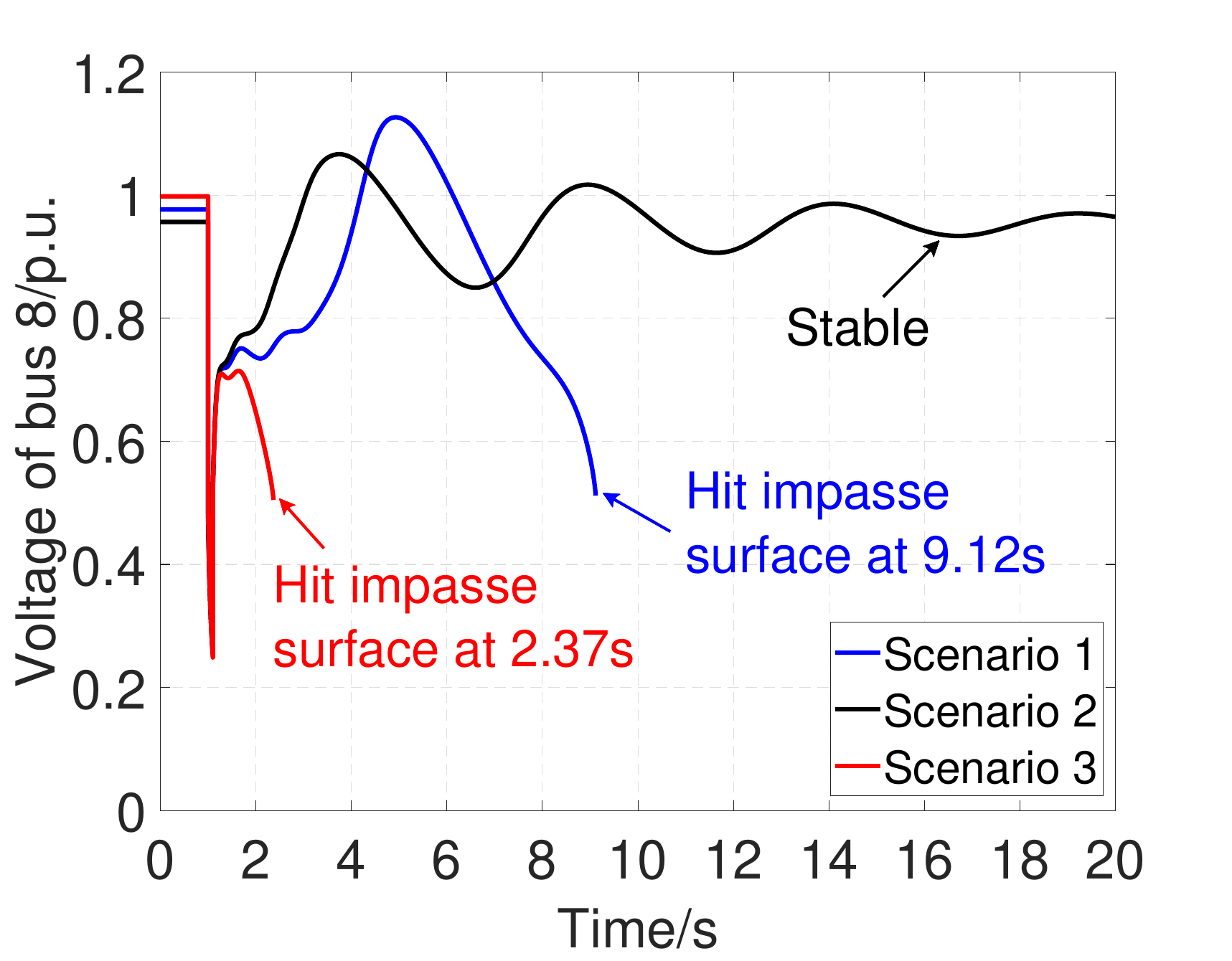}}
  %\hspace{0.1in}
  \subfigure[Min-modulus eigenvalue of $\Jalg$.]{
  \label{figcase9Jalg} %% label for second subfigure
  \includegraphics[width=1.67in]{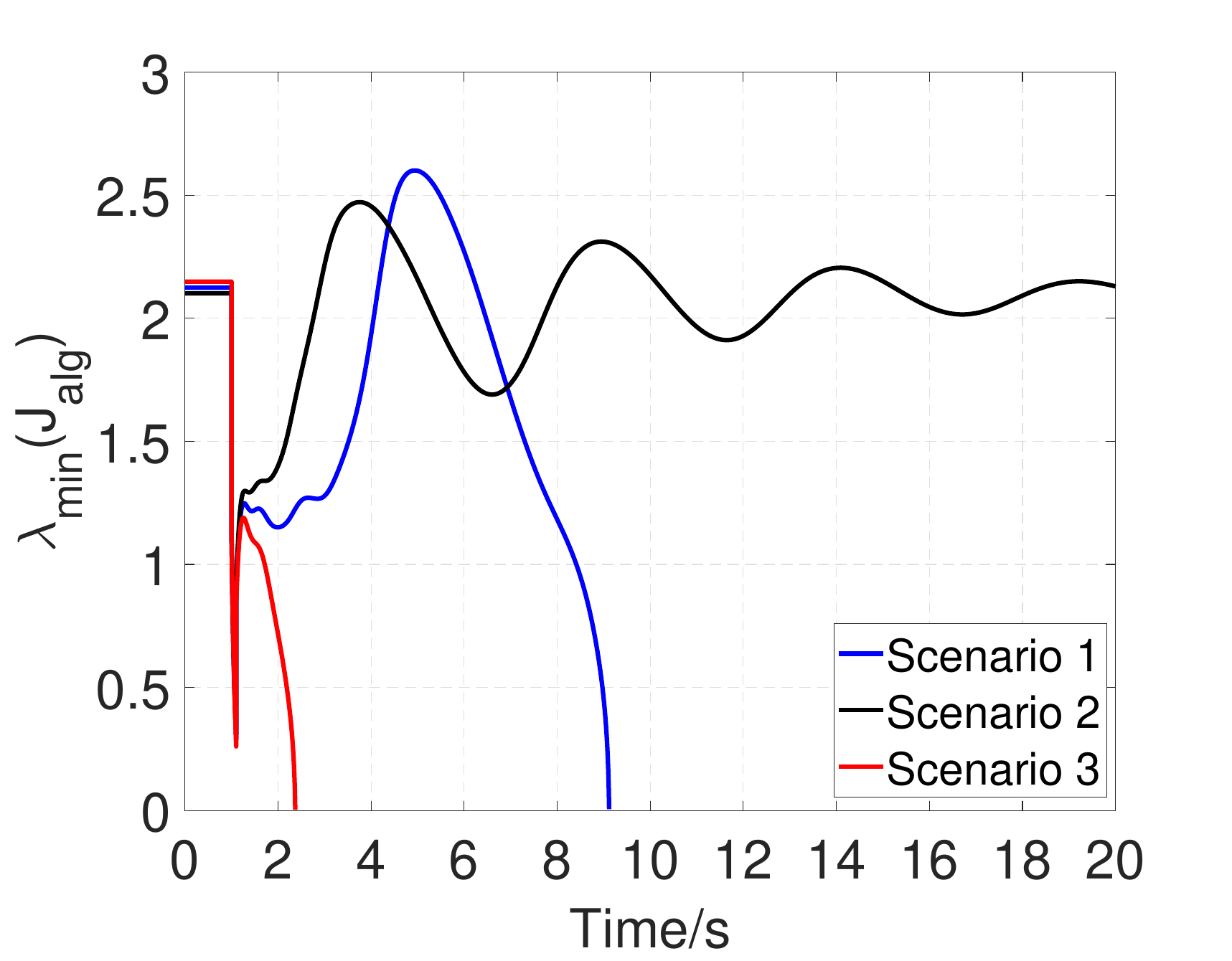}}
  \caption{System trajectories in three scenarios.}
  \label{figchargingschedule} %% label for entire figure
\end{figure}

\indent
We set the following three scenarios to verify the role of shunt capacitor/inductor in voltage collapse:
\begin{description}[leftmargin=1.8cm]
  \item[Scenario 1] The system parameters are as in \cite{case9data2019};
  \item[Scenario 2] A 0.30 p.u. shunt inductor is added to bus 8;
  \item[Scenario 3] A 0.30 p.u. shunt capacitor is added to bus 8.
\end{description}
For each of the scenarios, suppose the system initially operates at the stable equilibrium point and a three-phase short-circuit fault occurs at 1.0~s such that bus 8 is grounded via a 0.05~p.u. reactance, which is cleared at 1.1~s.
The system response with respect to bus 8 and minimum modulus eigenvalue of the algebraic Jacobian are depicted in Fig. \ref{figcase9voltage} and Fig. \ref{figcase9Jalg}, respectively.
Voltage collapse occurs in scenario 1 and scenario 3 at 9.12~s and 2.37~s, respectively, where the corresponding algebraic Jacobian becomes singular (see the blue and red curves in Fig. \ref{figcase9Jalg}). Contrarily, the post-fault system is stable in scenario 2.
It implies that the additional shunt inductor helps to prevent voltage collapse, while the additional shunt capacitor makes voltage collapse occur even earlier. This observation coincides with the analysis in Section \ref{secshunt}.

\begin{figure}[!h]
  \centering
  \subfigure[Scenario 1 and scenario 3.]{
  \label{figcase9Ivs13} %% label for first subfigure
  \includegraphics[width=1.67in]{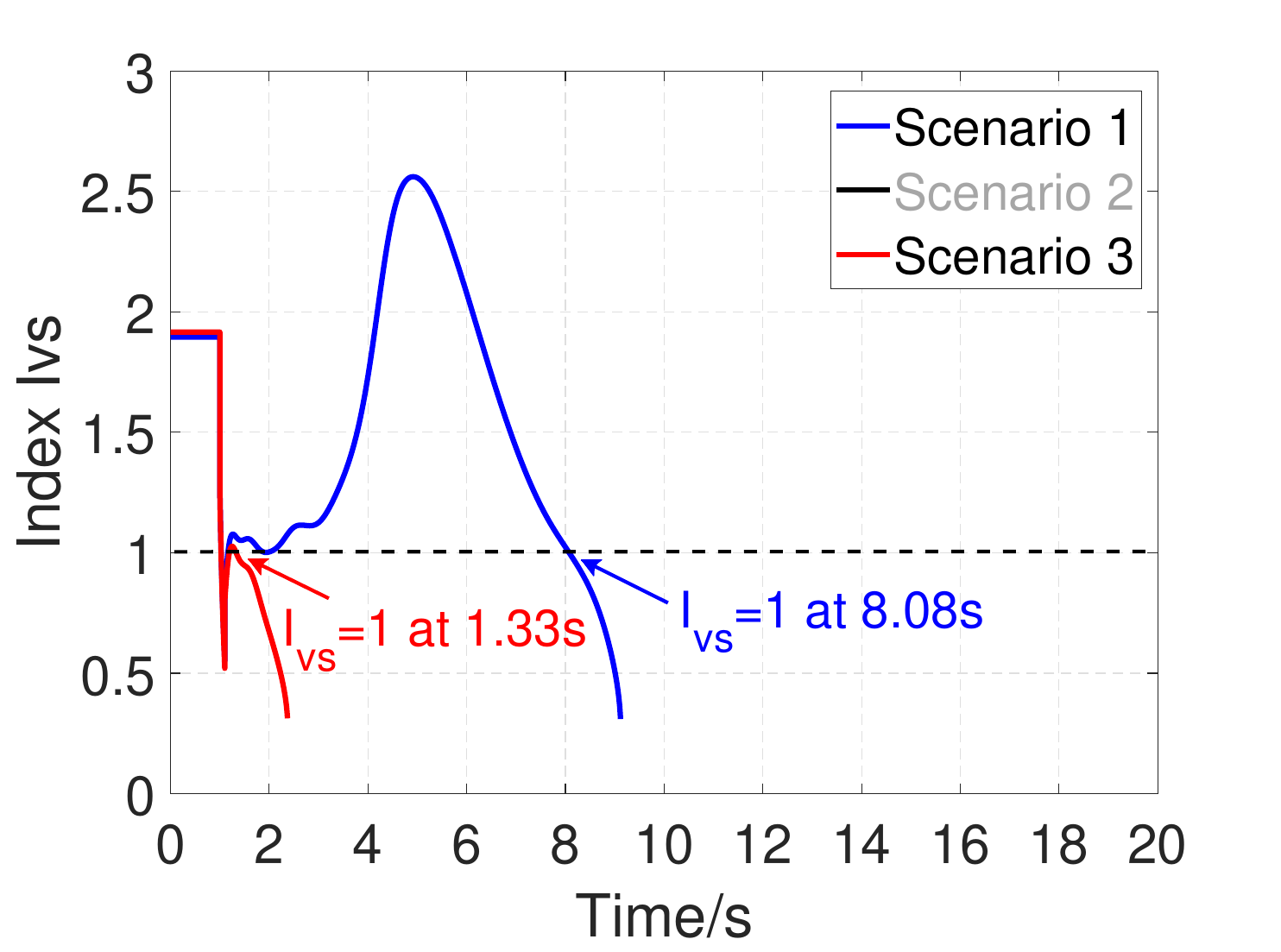}}
  %\hspace{0.1in}
  \subfigure[Scenario 2.]{
  \label{figcase9Ivs2} %% label for second subfigure
  \includegraphics[width=1.67in]{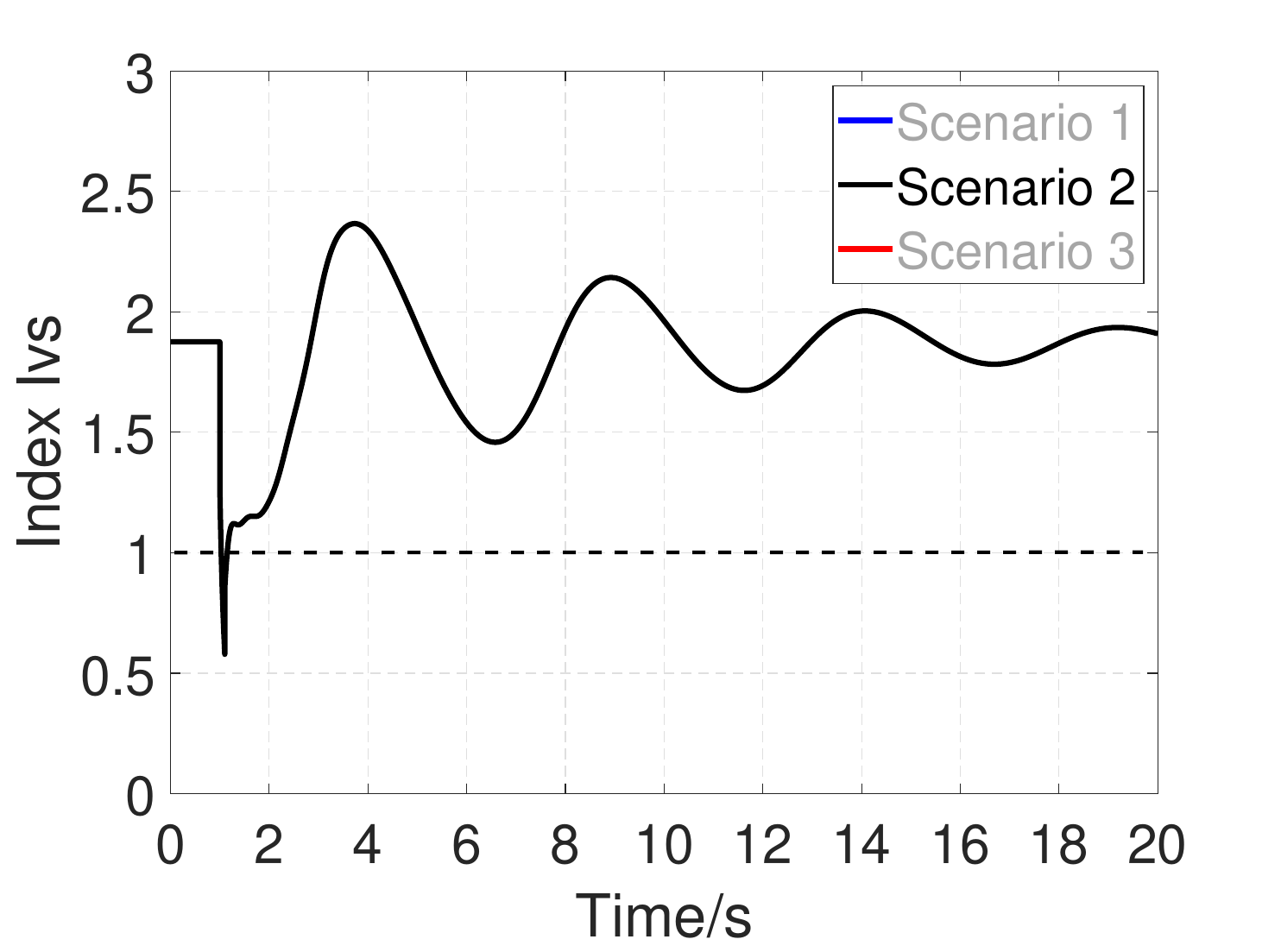}}
  \caption{Trajectories of $\Ivs(t)$ in three scenarios.}
  \label{figcase9Ivs} %% label for entire figure
\end{figure}

\indent
Now turn to the trajectories of the index $\Ivs(t)$ in Fig. \ref{figcase9Ivs} to illustrate inequality \eqref{impassenetwineq} in Theorem \ref{thmimpassenetw}.
For scenario 1 and scenario 3 that are unstable, $\Ivs(t)$ becomes greater than one for a short period due to temporary voltage recovery, but drops below one before voltage collapse (see Fig. \ref{figcase9Ivs13}).
For scenario 2 that is post-fault stable, $\Ivs(t)$ is less than one for a short period right after the fault is cleared, and converges to a steady-state value that is greater than one (see Fig. \ref{figcase9Ivs2}).
All these observations coincide with Theorem \ref{thmimpassenetw}.

\begin{figure}[!h]
  \centering
  \subfigure[Scenario 1.]{
  \label{figcase9shed808} %% label for first subfigure
  \includegraphics[width=1.67in]{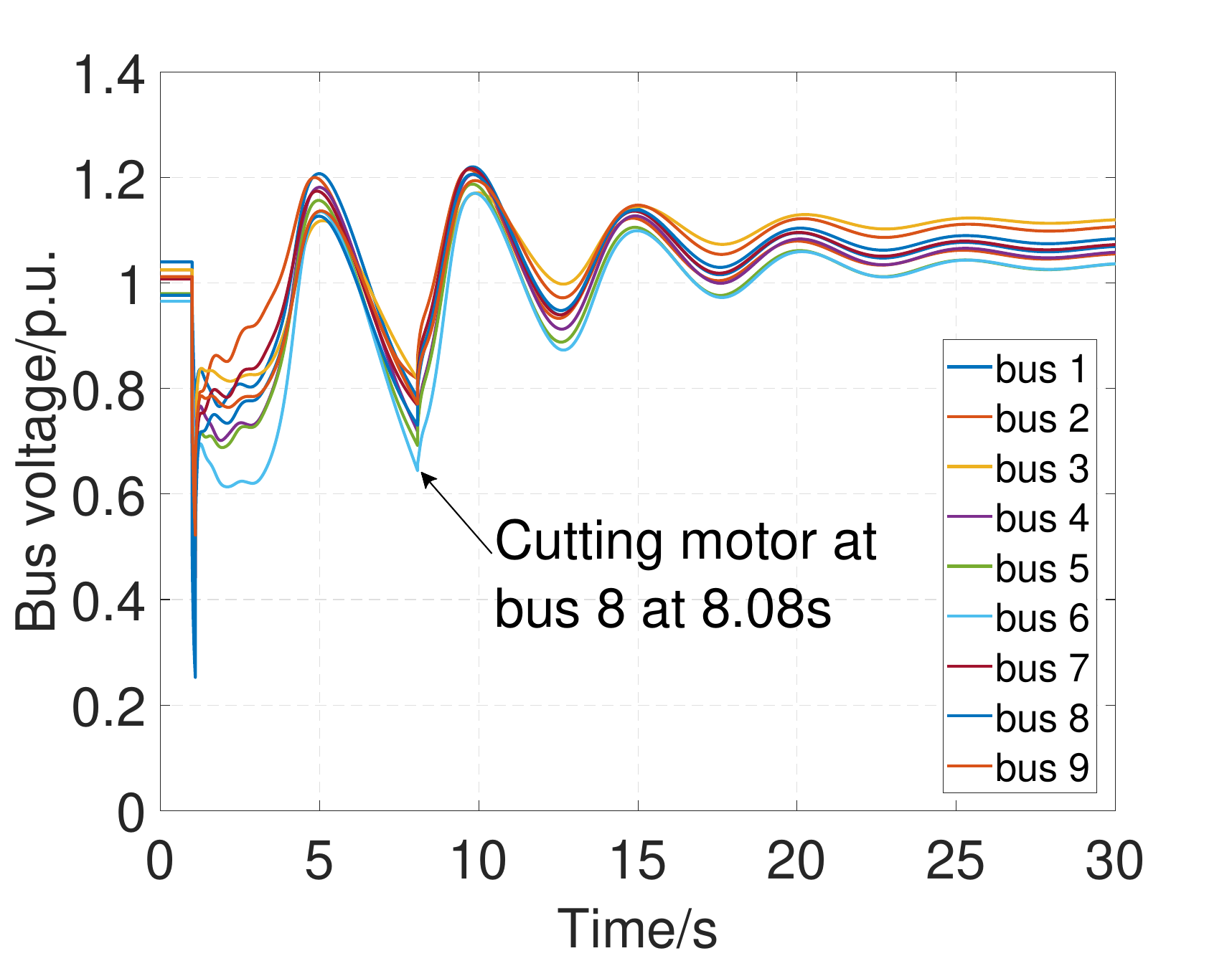}}
  %\hspace{0.1in}
  \subfigure[Scenario 3.]{
  \label{figcase9shed133} %% label for second subfigure
  \includegraphics[width=1.67in]{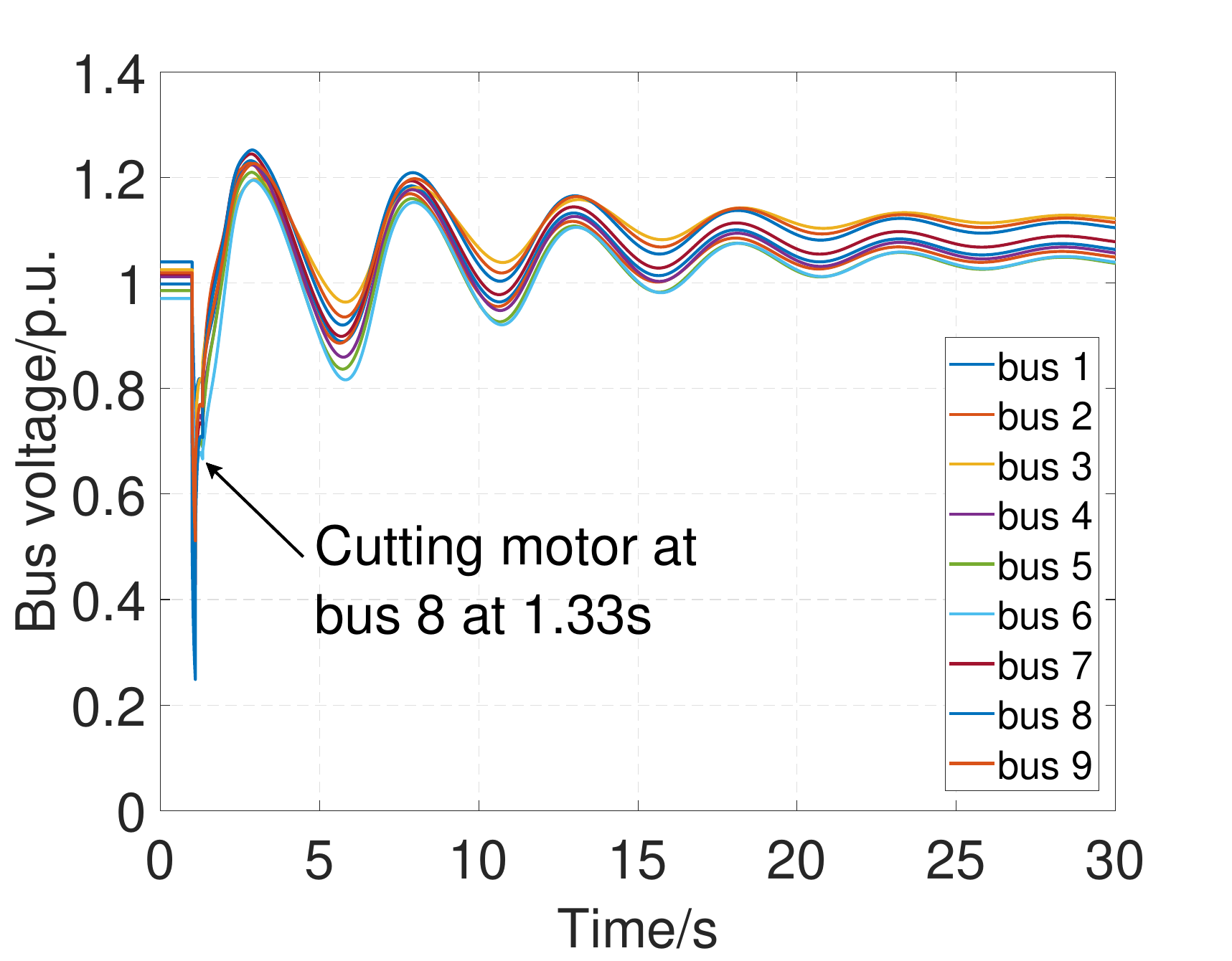}}
  \caption{Voltage trajectories with load shedding.}
  \label{figcase9shed} %% label for entire figure
\end{figure}

\indent
Moreover, we preliminarily show the potential of utilizing the index $\Ivs(t)$ in corrective control of voltage collapse.
As inferred from Theorem~\ref{thmimpassenetw}, ``$\Ivs(t)$ dropping below one'' serves as an early indicator of the system being closer to the impasse surface, which can be used to trigger control actions (e.g., load shedding) for preventing voltage collapse.
For instance, if the induction motor at bus 8 is cut when $\Ivs(t)$ drops to one, which respectively refers to 8.08~s and 1.33~s in scenario 1 and scenario 3 (the action time of cutting the motor is ignored for simplicity), then the corresponding post-fault system becomes stable (see Fig. \ref{figcase9shed808} and Fig. \ref{figcase9shed133}), which provides a new viewpoint for short-term voltage stability enhancement.
On the other hand, due to the conservativeness of Theorem~\ref{thmimpassenetw}, ``$\Ivs(t)$ dropping below one'' alone may not be an adequate criterion by itself in practice.
For instance, the stable trajectory in scenario 2 also experiences a very short period of time where $\Ivs(t)$ is below one (see the sharp sag around 1.0~s in Fig. \ref{figcase9Ivs2}).
It implies that some other logics, such as the duration of $\Ivs(t)$ being below one, need to be supplemented to achieve a better decision making on control actions.
A more systematic and practical framework on this issue will be further explored in the future.

\vspace{-4mm}
\vspace{-3mm}

\vspace{-2mm}

\vspace{-2mm}

\section{Conclusion}\label{secconclu}
A new characterization of the impasse surface of power system DAE model has been presented.
Admittance matrix-based necessary conditions for system trajectory hitting the impasse surface have been established, which reveal how the interactions between power network, generators and loads induce or prevent short-term voltage collapse.
The obtained theorems allow generic models for power network, generators and loads, which extend some existing results developed on simplified system models. In particular, our results prove the conjecture in \cite{hiskens1990phd} that had been pending for decades.
Moreover, the obtained theorems lead to an early indicator of voltage collapse and a novel viewpoint that inductive compensation has a positive effect on preventing short-term voltage collapse, which have been verified via numerical simulation on the IEEE 9-bus system.
Future works include a more comprehensive corrective control method based on these theorems and coordinating the reactive power compensation requirements for achieving both short-term and long-term voltage stability.
In addition, another conjecture in \cite{hiskens1990phd} says that the impasse surface is avoided if $P_{si}^0\geq0$, $\alpha_i\geq 1$, $Q_{si}^0\geq0$, $\beta_i=2$, $\forall i\in\mV_L$, which remains open and needs more studies.

\vspace{-3mm}
\section*{Appendix}
We first present the lemma below which serves as a basis for the proofs of Theorem \ref{thmimpassenetw} and Theorem \ref{thmimpasseexponent}.
\begin{lemma}\label{lemJalgYp}
     The algebraic Jacobian $\Jalg$ is singular if and only if the following matrix is singular
     \begin{equation}\label{Yp}
     \begin{split}
       \bm{Y}^{\prime}=
       \begin{bmatrix}
         \lbar{\bm{Y}_1} & \lbar{\bm{Y}_2\bm{T}} \\
         \bm{Y}_2\bm{T} & \bm{Y}_1 \\
       \end{bmatrix}\in\mbR^{2n\times 2n}
    \end{split}
    \end{equation}
     where $\bm{T}=\diag{e^{\upj 2\theta_i}}\in\mbR^{n\times n}$, $\forall i\in\mV_L$; $\bm{Y}_1$ is defined in \eqref{Y1}; $\lbar{\bm{Y}_1}$ denotes the entry-wise complex conjugate of $\bm{Y}_1$; and
     \begin{equation}\label{Y2}
     \begin{split}
       \bm{Y}_2(t) = (\bm{I}_n-\frac{1}{2}\bm{\alpha})\bm{G}_{\stat}^{\upeq}
       +\upj(\bm{I}_n-\frac{1}{2}\bm{\beta})\bm{B}_{\stat}^{\upeq}.
     \end{split}
    \end{equation}
\end{lemma}

\begin{IEEEproof}
   First, it is trivial that the singularity of $\Jalg$ is equivalent to that of the following matrix
   \begin{equation}\label{Jalgp}
   \begin{split}
      \Jalg^{\prime} =
   \begin{bmatrix}
     \frac{\partial \bm{g}_p}{\partial \bm{\theta}} & \frac{\partial \bm{g}_p}{\partial \bm{V}}\hat{\bm{V}} \\
     \frac{\partial \bm{g}_q}{\partial \bm{\theta}} & \frac{\partial \bm{g}_q}{\partial \bm{V}}\hat{\bm{V}} \\
   \end{bmatrix} \triangleq
   \begin{bmatrix}
     \bm{E} & \bm{F} \\
     \bm{M} & \bm{N} \\
   \end{bmatrix}
   \end{split}
   \end{equation}
   where $\hat{\bm{V}}=\diag{V_i}\in\mbR^{n\times n}$, $\forall i\in\mV_L$.

   \indent
   Let us further look into $\Jalg^{\prime}$.
   Observing \eqref{staticload}, \eqref{pfequ}, \eqref{Jpt} and \eqref{Ysi}, the entries of submatrices $\bm{E}, \bm{F}, \bm{M}, \bm{N}$ in \eqref{Jalgp} can be re-expressed in terms of $G_{mi}^{\upeq},B_{mi}^{\upeq}, G_{si}^{\upeq},B_{si}^{\upeq}$ as follows
   \begin{equation}\label{Jpt2}
   \begin{split}
       E_{ij} &=
       \left\{
        \begin{array}{ll}
           -V_i^2(\widetilde{B}_{ii}+B_{mi}^{\upeq}+B_{si}^{\upeq}),~i=j \\
           -V_iV_j|\widetilde{Y}_{ij}|\cos(\theta_{ij}-\varphi_{ij}),~i\neq j
        \end{array}
        \right.\\
     F_{ij} &=
     \left\{
        \begin{array}{ll}
           V_i^2(\widetilde{G}_{ii}+G_{mi}^{\upeq} + (\alpha_i-1)G_{si}^{\upeq}),~i=j \\
           V_iV_j|\widetilde{Y}_{ij}|\sin(\theta_{ij}-\varphi_{ij}),~i\neq j
        \end{array}
        \right.\\
     M_{ij} &=
      \left\{
        \begin{array}{ll}
           -V_i^2(\widetilde{G}_{ii}+G_{mi}^{\upeq} + G_{si}^{\upeq}),~i=j \\
           -V_iV_j|\widetilde{Y}_{ij}|\sin(\theta_{ij}-\varphi_{ij}),~i\neq j
        \end{array}
        \right.\\
     N_{ij} &=
       \left\{
        \begin{array}{ll}
           -V_i^2(\widetilde{B}_{ii}+B_{mi}^{\upeq} + (\beta_i-1)B_{si}^{\upeq}),~i=j \\
           -V_iV_j|\widetilde{Y}_{ij}|\cos(\theta_{ij}-\varphi_{ij}),~i\neq j.
        \end{array}
        \right.
   \end{split}
   \end{equation}
   We label the rows and columns of $\Jalg^{\prime}$ by the index set $\mI_0=\{1,2,...,n,1^{\prime},2^{\prime},...,n^{\prime}\}$.
   Let $\bm{E}_{r}\in \mbR^{2n\times 2n}$ be the elementary matrix that changes the row order from $\mI_0$ to a new one, say $\mI_1=\{1,1^{\prime},2,2^{\prime},...,n,n^{\prime}\}$.
   Then we obtain $\Jalg^{\prime\prime}$ by the following elementary transform
   \begin{equation}\label{Jalgpp}
   \begin{split}
      \Jalg^{\prime\prime}=\bm{E}_{r}\Jalg^{\prime}\bm{E}_{r}^{-1}.
   \end{split}
   \end{equation}
   By \eqref{Jpt2}, $\Jalg^{\prime\prime}$ can be expanded as
   \begin{equation}\label{Aij}
   \begin{split}
   \small
      \Jalg^{\prime\prime} =
      \begin{bmatrix}
        \bm{J}^{\prime\prime}_{11} &  \bm{J}^{\prime\prime}_{12} & \cdots  &  \bm{J}^{\prime\prime}_{1n}\\
        \bm{J}^{\prime\prime}_{21} &  \bm{J}^{\prime\prime}_{22} & \cdots  &  \bm{J}^{\prime\prime}_{2n}\\
        \vdots &  \vdots & \ddots  &  \vdots \\
        \bm{J}^{\prime\prime}_{n1} &  \bm{J}^{\prime\prime}_{n2} & \cdots  &  \bm{J}^{\prime\prime}_{nn}\\
      \end{bmatrix}
   \end{split}
   \end{equation}
   where $\bm{J}^{\prime\prime}_{ii}$, $i=1,2,...,n$ takes value as
   \begin{equation}
   \begin{split}
      \bm{J}^{\prime\prime}_{ii}=
      \begin{bmatrix}
        (\Jalg^{\prime})_{ii} &  (\Jalg^{\prime})_{ii^{\prime}} \\
        (\Jalg^{\prime})_{i^{\prime}i}  &  (\Jalg^{\prime})_{i^{\prime}i^{\prime}} \\
      \end{bmatrix}=
      \begin{bmatrix}
        E_{ii} &  F_{ii} \\
        M_{ii}  &  N_{ii} \\
      \end{bmatrix}
   \end{split}
   \end{equation}
   and $\bm{J}^{\prime\prime}_{ij}$, $i,j=1,2,...,n$, $i\neq j$, takes value as
   \begin{equation}
   \begin{split}
      \bm{J}^{\prime\prime}_{ij}=
      \begin{bmatrix}
        (\Jalg^{\prime})_{ij} &  (\Jalg^{\prime})_{ij^{\prime}} \\
        (\Jalg^{\prime})_{i^{\prime}j}  &  (\Jalg^{\prime})_{i^{\prime}j^{\prime}} \\
      \end{bmatrix}=
      \begin{bmatrix}
        E_{ij} &  F_{ij} \\
        M_{ij}  &  N_{ij} \\
      \end{bmatrix}.
   \end{split}
   \end{equation}

   \indent
   Let $\bm{U} = \frac{\sqrt{2}}{2}
      \begin{bmatrix}
        1 &  1 \\
        -\upj  &  \upj \\
      \end{bmatrix}$, we have
   \begin{equation}\label{Aij2Dij}
   \begin{split}
      \bm{U}^{-1}\bm{J}^{\prime\prime}_{ii}\bm{U} &= \upj
      \begin{bmatrix}
        -(\lbar{\bm{Y}}_1)_{ii}V_i^2 &  -(\lbar{\bm{Y}_2})_{ii}V_i^2 \\
        (\bm{Y}_2)_{ii}V_i^2  &  (\bm{Y}_1)_{ii}V_i^2 \\
      \end{bmatrix}\\
      \bm{U}^{-1}\bm{J}^{\prime\prime}_{ij}\bm{U} &= \upj
      \begin{bmatrix}
        -(\lbar{\bm{Y}}_1)_{ij}V_iV_j e^{\upj \theta_{ij}} & 0 \\
        0 &  (\bm{Y}_1)_{ij}V_iV_j e^{-\upj \theta_{ij}} \\
      \end{bmatrix}.
   \end{split}
   \end{equation}
   By \eqref{Aij} and \eqref{Aij2Dij}, $\Jalg^{\prime\prime}$ can be re-expressed as
   \begin{equation}\label{Jalgpp2D}
      \Jalg^{\prime\prime} = (\bm{I}_n\otimes\bm{U})\bm{K}(\bm{I}_n\otimes\bm{U})^{-1}
   \end{equation}
   where $\otimes$ denotes the Kronecker product and
   \begin{equation}\label{D}
   \small
      \bm{K} =
      \begin{bmatrix}
        \bm{K}_{11} &   \cdots  &  \bm{K}_{1n}\\
        \vdots  & \ddots  &  \vdots \\
        \bm{K}_{n1} &   \cdots  &  \bm{K}_{nn}\\
      \end{bmatrix}
   \end{equation}
   with $\bm{K}_{ij}=\bm{U}^{-1}\bm{J}^{\prime\prime}_{ij}\bm{U}$, $\forall i,j\in\mV_L$.
   We label the rows and columns of $\bm{K}$ by the index set $\mI_1$. Rearranging the rows and columns of $\bm{K}$ into the order $\mI_0$ gives the matrix $\bm{E}_{r}^{-1}\bm{K}\bm{E}_{r}$.
   Observing \eqref{Yp}, \eqref{Aij2Dij} and \eqref{D}, $\bm{E}_{r}^{-1}\bm{K}\bm{E}_{r}$ takes the form below
   \begin{equation}\label{Yp2D}
   \begin{split}
      \bm{E}_{r}^{-1}\bm{K}\bm{E}_{r}
      &=
      \begin{bmatrix}
         -\bm{I}_n & \bm{0} \\
         \bm{0} & \bm{I}_n \\
       \end{bmatrix}
      \begin{bmatrix}
        \hat{\bm{V}}^c & \bm{0} \\
        \bm{0} & \overline{\hat{\bm{V}}^c} \\
      \end{bmatrix}
      \bm{Y}^{\prime}
      \begin{bmatrix}
      \overline{\hat{\bm{V}}^c} & \bm{0} \\
        \bm{0} & \hat{\bm{V}}^c \\
      \end{bmatrix}
   \end{split}
   \end{equation}
   where $\hat{\bm{V}}^c=\diag{V_i e^{\upj\theta_i}}\in\mbC^{n\times n}$, $\forall i\in\mV_L$.
   From \eqref{Jalgp}, \eqref{Jalgpp}, \eqref{Jalgpp2D} and \eqref{Yp2D} we conclude that $\Jalg$ is singular if and only if $\bm{Y}^{\prime}$ is singular.
\end{IEEEproof}

\indent
Now we come to the proofs of the two theorems.

\indent
\textit{Proof of Theorem \ref{thmimpassenetw}: }
   Suppose \eqref{impassenetwineq} is violated, then we have
   \begin{equation}\label{blockdominant}
   \begin{split}
       \sigma_{\min}(\bm{Y}_1) > \sigma_{\max}(\bm{Y}_2\bm{T})
   \end{split}
   \end{equation}
   since the right-hand-side of \eqref{impassenetwineq} equals to $\sigma_{\max}(\bm{Y}_2\bm{T})$.
   By \eqref{blockdominant}, $\bm{Y}_1$ is non-singular, and hence \eqref{blockdominant} is equivalent to $\|\bm{Y}_1^{-1}\|_2^{-1} > \|\bm{Y}_2\bm{T}\|_2$, where $\|\cdot\|$ denotes the 2-norm of a matrix. It implies that $\bm{Y}^{\prime}$ is block strictly diagonally dominant \cite{feingold1962block} for the two-by-two block partition given in \eqref{Yp}.
   Thus, it follows from [\citen{feingold1962block}, Theorem 1] that $\bm{Y}^{\prime}$ is nonsingular.
   By Lemma \ref{lemJalgYp}, $\Jalg$ is also nonsingular so that the system trajectory will not hit the impasse surface.
\hfill $\blacksquare$

\indent
\textit{Proof of Theorem \ref{thmimpasseexponent}: }
     First we point out that $\bm{Y}_2=\upj(\bm{I}_n-\frac{1}{2}\bm{\beta})\bm{B}_{\stat}^{\upeq}$ under the given conditions.
     In addition, by the given conditions we have
     \begin{equation}\label{nodeineqinduct1}
     \begin{split}
         &|\upj B_{ii}+G_{gi}+\upj B_{gi}+\frac{1}{2}\alpha_iG_{si}^{\upeq}+\upj \frac{1}{2}\beta_iB_{si}^{\upeq}| \\
         &\geq |B_{ii}|+|B_{gi}|+|\frac{1}{2}\beta_iB_{si}^{\upeq}|\\
         &>\sum_{j=1,j\neq i}^n |B_{ij}|+|(1-\frac{1}{2}\beta_i)B_{si}^{\upeq}|,~\forall i\in \mV_t
     \end{split}
     \end{equation}
     which implies the rows of $\bm{Y}^{\prime}$ (defined in \eqref{Yp}) with respect to $\mV_t$ are strictly diagonally dominant.
     We also have
     \begin{equation}\label{nodeineqinduct2}
     \begin{split}
         &|\upj B_{ii}+\frac{1}{2}\alpha_iG_{si}^{\upeq}+\upj \frac{1}{2}\beta_iB_{si}^{\upeq}|\geq
         |B_{ii}| + |\frac{1}{2}\beta_iB_{si}^{\upeq}| \\
         &\geq\sum_{j=1,j\neq i}^n |B_{ij}|+|(1-\frac{1}{2}\beta_i)B_{si}^{\upeq}|,~\forall i\in \mV_L\backslash\mV_t
     \end{split}
     \end{equation}
     which implies the rows of $\bm{Y}^{\prime}$ with respect to $\mV_L\backslash\mV_t$ are diagonally dominant.

     \indent
     Further, we define a directed graph associated with $\bm{Y}^{\prime}$ as follows.
     The set of nodes of the graph is given by $\{1,2,...,2n\}$ and there is an edge orienting from node $i$ to $j$ if and only if $Y_{ij}^{\prime}\neq 0$.
     This directed graph is strongly connected as the physical power network (interpreted by $\bm{Y}_1$) is connected. %(see the example in Fig. \ref{figWCDD}).
     Thus, for any node $i\in \mV_L\backslash\mV_t$, there exists a path from node $i$ to $j$ in this directed graph such that $j\in\mV_t$.
     From the above discussion, $\bm{Y}^{\prime}$ satisfies the conditions to be weakly chained diagonally dominant (WCDD), and hence it is nonsingular \cite{azimzadeh2016weakly}.
     Then, by Lemma \ref{lemJalgYp}, any system trajectory will not hit the impasse surface.
\hfill $\blacksquare$

\ifCLASSOPTIONcaptionsoff
  \newpage
\fi

{\footnotesize
\bibliographystyle{IEEEtran}
\bibliography{IEEEabrv,voltageimpasse}

}


%\vfill


\end{document}